 \documentclass [11pt]{amsart}
\usepackage{amsmath}
\usepackage{xcolor}
\usepackage{amsfonts}
\usepackage{amssymb}
\newtheorem{theo}{Theorem} 
\newtheorem{prop}{Proposition}

\newtheorem{lema}{Lemma}
\newtheorem{defi}{Definition}

\begin{document}

\title{On the realization of the Gelfand Character of a finite  group as a twisted trace 
 }

%\author{Antonio Behn}
%\address{Dep. Mat. Fac. Ciencias, Univ. de Chile}

\author{Jorge Soto-Andrade}
\address{Dep. Mat. Fac. Ciencias, Univ. de Chile}

\author{ M. Francisca Y\'a\~nez } 
\address{Escuela de Pregrado, Fac. Cs. Qu\'imicas y Farmac\'euticas, Univ. de Chile}

\thanks{The authors have been partially supported by Fondecyt Grant  1140510}
\date{}

\begin{abstract}

We show   that the  Gelfand character $ \chi_G$ of  a finite group $G $ (i.e. the sum of all irreducible complex characters
   of  $G$ ) may be realized as a `` twisted trace''   $ g \mapsto   Tr( \rho_g \circ T)  $ for a suitable involutive linear automorphism $T$ of   $L^2(G)$, where   $(L^2(G), \rho)$  is the right regular representation of   $G$.  Moreover, we prove that under certain hypotheses   $T(f)= f \circ  L \;\; (f \in L^2(G)), $ where $ L $ is an involutive antiautomorphism
of $ G.$ The natural representation  $\tau$ of $G$ associated to the natural $L$-conjugacy action of $G$ in the fixed point set  $Fix_G(L)$ of   $L$   turns out to be a Gelfand Model for $G$ in some cases. We show that  $(L^2(Fix_G(L)), \tau) $   fails to be a  Gelfand Model  if    $G$ admits non trivial central involutions.

\end{abstract}

 \maketitle

\smallskip

\smallskip


\section{Introduction}

%%% Check notations and recall isotypical  components 

In the seventies, Gelfand and collaborators \cite{bgg} introduced the notion of a {\em Model} for the representations of a  compact Lie group $G$, as a representation of $G$ where each irreducible representation of $G$ appears exactly once.  Nowadays these Models are usually callled {\em Gelfand Models}.  

A paradigmatic example of such a Gelfand Model is the natural representation of  the compact Lie group $G  = SO(3, \mathbb R)$ in the Hilbert space $  H = L^2(S^2),$ where $S^2$ denotes the unit sphere in $\mathbb R^3.$ The spherical functions associated to the (Hilbert sum) decomposition of this representation into irreducible mutually orthogonal subspaces are given by the well known Legendre Polynomials $P_{\ell} \;\;  (\ell \in \mathbb N).$ More precisely, the orthogonal projectors onto the irreducible components  $ H_\ell $ of  $H $ are integral operators with kernels  $K_\ell (x,y) = P_\ell (x \cdot y) \;\; (x,y \in S^2)$.

This suggests that remarkable special functions will appear as spherical functions providing  the decomposition into irreducibles of Gelfand Models for particular groups.

The notion of Gelfand Model makes sense of course for any finite group $G,$ in which case a Gelfand Model is a representation of $G$ isomorphic to the direct sum of all  (complex) irreducible representations of $G.$ 

It is however very rarely the case (as it is for $SO(3, \mathbb R)$) that a finite {\em non commutative} group admits a Gelfand Model which is just a natural representation associated to some  (necessarily transitive) $G-$ space $X.$   In that case we   say that $G$ admits a {\em natural Gelfand Model}.  

Nevertheless, as foreseen by Gelfand, these representation have often remarkable properties. 
 Indeed, this the case for example for $G = PGL(2, \mathbb F_q),$ for which the main constituent of a Gelfand Model  is the natural representation of $G$ in 
 $ L^2(  \widetilde {\mathcal H}) $         
associated to the homographic action of $G$ on  $ \widetilde {\mathcal H} =  \mathbb F_{q^2} -  \mathbb F_q, $  the double cover of the finite analogue of Poincar\'e's Upper Half Plane over the finite base field $ \mathbb F_q$ \cite{saarc, jajatwp2}. The spherical functions providing the decomposition of    $ L^2(  \widetilde {\mathcal H}) $  into irreducibles  have turned out to play a crucial role in the construction  of   regular Ramanujan graphs associated to  $ \widetilde {\mathcal H} $ \cite{katz, winniebook}.
  
  Notice that we keep  the $ L^2$ notation in the finite case, so that $L^2(X) $ for a finite set $X$ denotes the complex vector space of all complex functions on $X,$  endowed with the usual   $L^2$ product.    
 
 In this paper we consider the character of a Gelfand Model of a finite group $G,$ which we call the {\em Gelfand Character} of $G$  and denote by $\chi_G.$   This character, which is just the sum of all irreducible character of $G$, has indeed remarkable properties.   
 
 It looks very much like a (transitive) permutation character (the character of a natural, or permutation, representation associated to a transitive $G$-space): Its mean value is 1 and  it is integer valued. Moreover, its values are very often non negative. In fact this is the case for the the ordinary finite classical   groups and we conjectured it  to be true for every finite group until Yokonuma \cite{yoko} pointed out to us that  the Gelfand character of the exceptional Mathieu group $M_{11}$  took negative values at some conjugacy classes. Later, Behn \cite{absa} found a minimal counter-example for this conjecture, to wit a group $G(96)$ of order 96 whose Gelfand Character takes values  [30, 2, -2, 6, 6, 2, 2, 2]  on its 8 conjugacy classes \cite{absa}.    
 
 Also, we may notice that  the alternating group $A_4$ admits a natural Gelfand Model  $ M = L^2(X),$ where the $G-$ set $X$ stands for the set of the six oriented edges of the regular tetrahedron, endowed with the natural action of $A_4.$  An interesting open question in this regard is to characterise all finite groups admitting a natural Gelfand Model.    
  
 We remark that a first obstruction for this to happen is that the dimension $d_G$ of the Gelfand Model of $G$ must divide the order of $G$, which is rarely the case.  However, remarkably enough, for $G = SL(2, \mathbb F_q)$ we have that the dimension of its Gelfand Model is $q(q+1)$, which divides the order $(q-1)q(q+1)$ of $G$, although as we show below there is no natural Gelfand Model for $G.$
 
The realization of the Gelfand character    $\chi _G$ of  a finite group  $G$, is an old problem \cite{bg, saarc, terras, mfycras}.  One possible approach to this problem is to try to obtain   $\chi _G$ by twisting the  trace of some very natural representation $(V, \pi) $ of   $G$, like its regular representation, by a suitable linear automorphism  $T$  of its underlying space $V,$ so as to obtain    $  \chi_G(g) = Tr(\pi _g \circ T) $   for all  $ g \in G.$
Recall that twisted traces appear in many contexts in mathematics  \cite{arhis,bg, ha}. 

A classical result in character theory  of finite groups \cite{gow} is that the central function   $ \theta_1: G \to \mathbb{C},$ defined by 
 \mbox{$ \theta _1(g)= |\{h \in G : h^2 = g \}| $} is a generalized character that  satisfies $ \theta _1(g)= \sum _{ \pi \in \hat{G}} \nu ( \pi ) \chi _{\pi}(g) $ where 
 $ \nu (\pi) = \frac{1}{|G|} \sum_{g \in G } \chi_{\pi}(g^2) $ is the Frobenius- Schur  number of the character $ \chi_{\pi} $  of the irreducible representation $ \pi $ of $ G. $
 If $ (\pi, V ) $ is an irreducible complex representations of a finite group $ G $ and $ \pi $ is self-contragradient  then there exists a
non degenerate  bilinear form $ B $ on $ V, $ unique up to scalar multiple, such that $ B(v,w) = \epsilon _(\pi) B(w,v) $ where $  \epsilon _(\pi) = \pm1. $ Frobenius  and Schur \cite{f-s} proved that
$ \nu (\pi) = \epsilon (\pi) $ and if $ \nu (\pi) = 1 $ then $ \pi (G) $ is conjugated to a subgroup of the orthogonal group $ O(n) $ and if $ \nu(\pi) = -1 $  then $ \pi (G) $ is conjugate to a subgroup of the symplectic group $ Sp(n) $ and  $n $ is even.
 Clearly if $ \nu(\pi)= 1 $ for all  irreducible representations of $ G,$ then $ \theta_1$ is the Gel'fand character $\chi_G$ of $G$. Some groups for which this is the case are $ S_{n}, D_{2n} $ and $ O(n, q ),  \;\;   q $ odd \cite{gow}. 
 
Furthermore  Gow \cite{gow} proved  that  for $ G= GL(n,q), $ the central function $ \theta_2: G \to \mathbb{C},$ defined by $ \theta_2 (g)= |\{h \in G : h^{t} h^{-1} = g \}| $ affords $\chi_G.$

On the other hand,  Kawanaka and Matsuyama  \cite{ K-M} defined the twisted Frobenius Schur indicator $ \nu _{\tau}( \chi ) =  \frac{1}{| G| }\sum _{g \in G } \chi ( g \tau (g))$ where $ \tau $ is an involutive automorphism of $G. $ 
 They  proved that if  the character $ \chi $ is afforded by a matrix representation $ R $ of $G $  such that  $ R \circ \tau = \overline{R }$   then $ \nu _{\tau}( \chi ) = 1.$  If this is not the case but  nevertheless  $ \chi \circ \tau  =  \overline {\chi },  $ then  $ \nu _{\tau}( \chi ) = -1.$ Finally if    $ \chi \circ \tau  \neq  \overline{\chi  }, $ then  $ \nu _{\tau}( \chi ) = 0.$ 
Furthermore they get that  $ \sum_{ \chi \in \hat{G}} \nu _{\tau } ( \chi ) \chi (g)  =  | \{ h \in G : \tau(h)^2 =g \}|.$ 

Later  Bump and Ginzburg \cite{bg} considered an automorphism $ \tau$ of  $ G $ such that $ \tau ^{r} = 1,  $ the norm map $ N: G \to G $ given by 
$ N(g)= g \tau (g)\tau ^2(g) \cdots \tau ^{r-1}(g) $ and  the number $ M(g) $ of solutions of the equation $ N( x ) = g $ with $ x \in G.  $ They proved that $ M(g) $ lies in the 
$ \mathbb{Z} e ^{\frac{2 \pi i }{r}} $- algebra generated by the irreducible characters and $ M(g) = \sum _{ \chi \in \hat{G} } \epsilon ( \chi ) \overline{\chi (g)}, $ where
$ \epsilon_{ \tau} ( \chi ) = \frac{1}{| G |} \sum _ { g \in G } \chi (N(g)). $ 
If $ r= 2 $ and $ \tau $ is trivial then $ \epsilon _ { \tau }( \chi ) = \nu ( \chi ) $ and $ M(g) = \theta _1(g) $; if   $ \tau$ is not trivial then $ \epsilon_{\tau}( \chi ) = \nu _ {\tau} ( \chi ), $ and $ M(g) = \sum _ { \chi \in \widehat{G} } \nu _{\tau}( \chi ) \chi (g). $ If $ G = GL(n, q) $ and $ \tau $ is  the transpose inverse then $ M(g) = \theta _2(g) .$

We prove that the Gelfand character $ \chi _{G} $ of a finite group $G$ may be always realized as a twisted trace  as $ \chi _{G}(g)  =  Tr( \rho _{g} \circ T ), \;\; g \in G, $ where $ T $ is an involutive automorphism  of $ L^2(G) $ and $  (L^2(G), \rho) $ is the right regular representation. 
Moreover we prove that if $G$ admits  an involutive anti-automorphism  $ L$  such that $ L(g) $ is conjugated to $ g, $  for all $g \in G$  then, putting     $   L^{*} (f) = f \circ L ,$   we get that  the central function 
 $ t(g) =   Tr (\rho _g \circ L^{*}) $ of    $g \in G $  equals the number of solutions $h \in G$ to the equation $ h^{-1}L(h) = g, $  and  is given by 
\mbox{$   \sum _{  \pi \in \widehat{G}} \varepsilon _{\pi} \chi _{\pi}(g), $}
  where 
$$  \varepsilon_{\pi }= \nu_{L}(\chi _{\pi})= \frac{1}{|G|}\sum_{ g \in G } \chi _{\pi} (L(g) g^{-1}) = \pm 1.$$

Furthermore if the number of fixed points of $ L $ is equal to the sum of the dimensions of all irreducible representations of $ G, $ then under the above assumptions the central function $ t $ is the Gelfand character of $G. $  So we can take the involutive automorphism $T$  of   $L^2(G) $ above to be  $ L^{*}.$

Let $ G $  be a group for which there exists  an anti-automorphism $L$  satisfying all the  conditions above.
Then we can consider the  $ G $ - set  $   X  =   \{ h \in G  :  L(h) = h      \} $  where the action of $G$ in $X$ is given by      
  $ g \cdot h = g^{-1}hL(g)^{-1}. $

Let   $\tau$ be the associated natural representation of  $G$ in $L^2(X)$.  Although  $ L^2(X)$ has the right  dimension to be a  Gelfand Model of  $G$,  we prove - using our twisted trace description of   $\chi_G$ - that this cannot be the case if  $G$ admits non trivial central involutions.

\section{The Gelfand character $\chi _G $ as a Twisted Trace}

Let  $G$ be a finite group and let  $ ( L^2(G),\rho)$ and $( L^2(G), \sigma )$ be the right and left
regular representation of $ G $ respectively;  let  $ (U_k , \pi _k ) \   ( 1
\leq k \leq r ) $ denote the irreducible unitary representations of $ G $ (up to isomorphism), with 
$n_k = {\rm dim} \; U_k$ and   let $ \emph{\textbf{I}}_k  \  ( 1 \leq k \leq r ) $ be the isotypic component of type $
\pi _k $ of $ \rho. $  We denote $ \chi _k $    the character of
 $ \pi _k .$

Let   $ ( e_{ij}^k(g)) _{1 \leq i,j \leq n_k} $ be the matrix of the operator $
\pi _k (g) \ ( g \in G) $ with respect to  an orthonormal basis
%${\mathcal U} _k, $ 
of $   U _k. $
%by $ \pi _k .$
%$ {\mathcal U} _k $ be an orthonormal  basis of $ U_k $ and the function $e_{ij}^k \in L^2(G)$ %defined by$$e_{ij}^k(g)=\langle \pi_k (g) (u_i), u_j\rangle , g \in G,  1 \leq i,j \leq n_k,  1
%\leq k \leq r $$
Then the
matrix coefficients $ e_{ij}^k  \ ( 1 \leq i,j \leq n_k, 1 \leq k \leq r)$ provide  an
orthonormal basis $ \mathcal{B} $ for the Hilbert space  $ L^2(G). $. More precisely, the matrix coefficients $ e_{ij}^k  \ ( 1 \leq i,j \leq n_k)$ 
provide an orthonormal basis for the isotypical component $I _k, $ for $ 1
\leq k \leq r  $  and for each  $1 \leq i \leq n_k $ the coefficients $  e_{ij}^k  \ ( 1 \leq j \leq n_k)$ afford an isomorphic copy of the irreducible representation  $ (U_k, \pi_k) $ of $G$ in $L^2(G).$

 The matrix coefficients $ e_{ij}^k $  satisfy the relations
 
\begin{equation}
e_{ij}^k(g^{-1})= \overline{e_{ji}^k(g)}
 \end{equation}

 and
 
\begin{equation}
e_{ij}^{k}(gh) = \sum_{ l = 1}^{ n_{k}}e_{il}^k(g)e_{lj}^k(h)
 \end{equation}
 
 for all $ 1 \leq i,j \leq n_k, 1 \leq k \leq r.$
 Moreover,  $ \chi _k $ denotes  the character afforded
by $ \pi _k .$
\begin{prop}
Let $ (V,\pi^1) $ and $( V,\pi^2) $ be two isomorphic representations
of a finite group $ G $ in the same vector space $V$ such that $ \pi^1_h \circ \pi^2_g = \pi^2_g
\circ \pi^1_h $ for all $g,h \in G $ and let $ S $ be an involutive automorphism  of $ V $
that intertwines the representations $ \pi^1 $ and $\pi^2 $ of $
G.$
Then the  complex function $ Tr( \pi^1_? \circ S ) = Tr ( \pi^2_? \circ S )$    defined on $ G $
is  central and so a linear complex combination of
irreducible characters of $ G. $
\end{prop}
\begin{proof}

For $ g, h $ in $ G, $ we have 

$$ Tr( \pi^1_{g^{-1}hg} \circ S ) = Tr( \pi^1_g \circ S \circ \pi^1_{g^{-1}h}) =   Tr(  S \circ \pi^2_g \circ \pi^1_{g^{-1}h}) =   $$
 
$$ = Tr( S \circ \pi^1_{g^{-1}h} \circ  \pi^2_g ) = Tr ( \pi^2_g \circ S \circ \pi^1_{g^{-1}h} )=Tr ( S \circ \pi^1_g   \circ \pi^1_{g^{-1}h} ) = 
$$ 

  $$ =Tr (S \circ \pi^1_h )   =Tr (\pi^1_h \circ S ).  $$

\end{proof}

 Recall that all irreducible representations $ \pi_k $ of a
finite group $ G $ are  unitarizable.  Moreover, the
isotypic component $ \emph{\textbf{I}}_{\pi_
k}$ of type $ (U_k,\pi_k) $ of $ \rho $ is isomorphic to $ U_k \otimes U_k^\ast. $

 \begin{theo}
 Let $ T $ be the  linear application of $ L^2(G) $ defined by $ T(
 e_{ij}^k) = e_{ji}^k  $ for all $ e_{ij}^{k} \in \mathcal{B}$ and let  $\widetilde{\sigma}$ be the homomorphism
from $ G $ to $Aut(L^2(G)) $ defined by
  %% gggggg 
  $$\widetilde{\sigma}_{g} = T \circ \rho_g \circ T $$ for all $ g \in G.$ 
  Then $ T $ is a linear   involutive automorphism
  of $  L^2(G) $ and
$\widetilde{\sigma}$ is a  representation  of 
 $ G $ such that:

\begin{itemize}

 \item[i.]      \ \ $ \rho _g \circ T = T \circ \widetilde{\sigma} _g, \;\;g \in G. $

 \item[ii.]
\ \ $ \rho _g \circ \widetilde{\sigma} _h = \widetilde{\sigma} _h \circ \rho _g,  \;\; g, h \in G. $
% ****.   ES BIEN SABIDO QUE  rho y sigma conmutan    
\item[iii.]
 \ \ $ Tr( \rho_g \circ T ) = \chi_{G} (g),$ \;\;
$g \in G $
 \end{itemize}
\end{theo}

\begin{proof}

  Since  $T$  is involutive, we get from \;\;
  $ \widetilde{\sigma}_{g} = T \circ \rho_g \circ T, \;$   
  that $\widetilde{\sigma}_g $ is an
  automorphism of $ L^2(G) $ and 
  that $ \rho_g \circ T = T \circ \widetilde{\sigma_g},$ for all  $ g \in G. $ 
  
  Furthermore, since 
  
   $$ {\rho}_{g}(e_{ij}^k) = \sum_{l=1}^{n_{k}} e_{lj}^k(g) e_{il}^k $$
   and 
   
  $$\widetilde{\sigma}_{g}(e_{ij}^k) = \sum_{l=1}^{n_{k}}e_{li}^k(g)e_{lj}^k $$ 
  for all $  e_{ij}^k \in \mathcal B ,    g \in G,$  we get that 
  
  %\noindent  

  $$ (\rho_g \circ \widetilde{\sigma}_h)(e_{ij}^k) =
  \sum_{ l = 1}^{ n_{k}}e_{li}^{k}(h)\sum_{m=1}^{n_{k}}e_{mj}^{k}(g)e_{lm}^k $$
   
  $$=\sum_{ m=1}^ { n_{k}}e_{mj}^{k}(g)\sum_{ l=1}^ { n_{k}}e_{li}^{k}(h)e_{lm}^k$$
  $$=\sum_{ m=1}^ {n_{k}}e_{mj}(g)\widetilde{\sigma}_{h}(e_{im}^k) $$ 
  
  $$=  (\widetilde{\sigma}_{h} \circ \rho_g)(e_{ij}^k).$$
 for $ g, h \in G $ and $ e_{ij}^k \in \mathcal{B}, $ which proves ii.
Finally, since 
  $$ (\rho_g \circ T)(e_{ij}^k)
  =\sum_{ l=1}^ { n_{k}}
  e_{li}^{k}(g)e_{jl}^k $$ 
  for all $ g \in G$ and $ e_{ij}^k \in \mathcal{B},$ then 
  $$ Tr(\rho_g \circ T) =
 \sum_{ k=1}^{r}(\sum_{ i =1}^{ n_{k}}e_{ii}^{k}(g) =
  \sum_{k=1}^{ r}\chi
 _{k}(g) =
  \chi_{G} (g),$$
  which proves iii.
\end{proof}

Next we will prove that under certain conditions the central function
$ Tr(\rho_{\it?} \circ T )$ can constructed from  an involutive anti-automorphism
  $ L $ of the group $ G. $
  
  % Notice that  counterexample for the existence of an L with the properties below, is provided %by the quaternion group H which has 5 irred repns.  (see Antonio)
  
  %Check that if. T is a lift of an anti autom L then L must behave prima ! 
  
  %Relocate remark 2 at the end.  
\begin{theo}

Let $ L $ be an involutive antiautomorphism   of $G,$ such that $
L(g) $ is conjugated to  $ g,$ for all $ g \in G, $ and let $ L^\ast $ be the
automorphism of $ L^2(G) $ defined by $ L^\ast (f)=f \circ L. $    Then for all $
g \in G, $ we have

\begin{itemize}
\item[i.]
$\rho_g \circ L^\ast  =   L^\ast  \circ \sigma_{L(g)^{-1}} $

 \item[ii.]
$ Tr(\rho_g \circ L^\ast ) = |\{h \in G : h^{-1}L(h)) = g \}|  $

 \item[iii.]
$ Tr(\rho_g \circ L^\ast ) =
 Tr(\rho_{g^{-1}} \circ L^\ast ) $
 \item[iv.]
 $ Tr(\rho_g \circ L^\ast ) = \Sigma_{k=1}^ { r} \varepsilon_k \chi_k (g) $
where $ \varepsilon_k = \pm 1. $
  .
\end{itemize}

\end{theo}

\begin{proof}
% $\rho_g \circ \sigma_{L(g)^{-1}} = \sigma_{L(g)^{-1}} \circ \rho_g$

The proof of i. is a straightforward
 calculation. We obtain ii. by computing the trace of \; $
 \rho_g \circ L^\ast  $ with
 respect to the canonical basis $ \{\delta_g : g\in G \} $ where $
 \delta_g (h) = \delta_{g,h}\; (h \in G)   .$ 

Now, let $ A_g =\{h \in G :  h^{-1}L(h) = g \}. $ Taking inverses we see that   $    h \in A_g $  if and only if $  L(h) \in A _{g^{-1}}. $
Then iii. follows since  $ L $ is bijective.    

To prove iv. notice first that because of  i.  the representations $  (L^2(G), \rho )$ and
 $ (L^2(G),  \check {\sigma}  ), $ where $ \check {\sigma}  (h) = \sigma_{L(h)^{-1}} \;\; (h\in G),   $
 together with  the involutive automorphism $L^{\ast} $ 
  satisfy all assumptions of  proposition 1, and so it follows 
  that the complex function $Tr(\rho_g \circ L^\ast ) $  is central and  
 $$ Tr(\rho_g \circ L^\ast )   = \sum_{ k=1}^ {r}\lambda_k \chi_k (g), $$ for suitable complex numbers 
 $\lambda_k.$

 The antiautomorphism $ L $ induces an antiautomorphism ${\tilde{L}}$
 on the complex group algebra $ \mathbb{C}[G]. $ Since $
\chi_{k}(L(g)) = \chi_{k}(g), g \in G, 1 \leq k \leq r,$  $
{\tilde{L}} $ acts as the identity on the centre  and therefore
induces an antiautomorphism $ \tilde{L}_k $ on each simple
component of $ \mathbb{C}[G] \cong \bigoplus _{ 1 \leq k \leq r }M(n_k ,\mathbb{C}).$
 Due to Skolem- Noether theorem, $ \tilde{L}_k $ is conjugated to the
 transposition map, i. e. there exists $b_k  \in GL(n_k,\mathbb{C})$ such that  $\tilde{L}_k(a) = b_k^{-1}a^{t}b_k $ for all $ a \in M(n_k,\mathbb{C}) $ . Furthermore 
we have that $ a = {\tilde{L}}({\tilde{L}}(a)) =
 b_k^{-1}b_k^ta({b_k^t})^{-1}b_k = b_k^{-1}b_k^t a(b_k^{-1}b_k^t)^{-1}$ for all $a \in M(n_k,\mathbb{C}),$ so $ b_k^{-1}b_k^t $ belongs to the centre  and
 therefore $ b_k^t = \varepsilon_k b_k $ with $\varepsilon_k = \pm 1,$ since $ ((b_k)^t)^t = b_k.$

 In this way, for each representation $(U_k,\pi_k) $ of $ G $
 a  symmetric     or a antisymmetric  form $ b_k $ exists, with respect to which
 the linear operators $ \pi_k(g) $ and $ {\tilde{L}}(\pi_k(g)) $     are adjoint to each other. More precisely, if we consider the  bilinear form $\langle u,v\rangle = v^t b_k u,$  defined by $b_k,$ then 
\begin{equation}
 \langle \tilde{L} (\pi_k(g))(u),v\rangle= v^t (b_k\tilde{L} (\pi_k(g)))u= v^t (\pi_k(g))^t bu=\langle u, \pi_k(g)(v)\rangle 
\end{equation}

Assume that for some  $k \; (1 \leq  k \leq r)$  we have $ \varepsilon_k = 1. $    We 
choose then an orthonormal basis $\mathcal{U}_k^{+}= \{ u_i : 1 \leq i
\leq n_k \} $ of $ U_k $, with respect to the symmetric form $ b_k $ and we denote
by $ e_{ij}^k(g)=\langle \pi_k (g) (u_j), u_i \rangle $  the matrix
coefficients of $ \pi_k(g) $ with respect to this basis. 
Due to equation 2 we have the following relations between  the matrix coefficients of 
$\tilde{L} (\pi_k(g)) = \pi_k(L(g))$ and $ \pi_k(g),$

\begin{equation}
   e_{ij}^k(L(g))  = \langle u_j,\pi_k(g)(u_i)\rangle = e_{ji}^k(g).
\end{equation}

Let $E^k= \langle e_{ij}^k: 1\leq i,j \leq n_k\rangle $ and $Tr_k(\rho_g \circ L^\ast )$ denote the trace of the restriction of $\rho_g \circ L^\ast $ to the subspace $E^k$  of $L^2(G).$  
In order to compute $ Tr_k(\rho_g \circ L^\ast ) $ we note that

$$ [( \rho_g \circ L^\ast )(e_{ij}^k)] (h)=e_{ij}^k(L(hg))=e_{ij}^k(L(g)L(h)) = $$
 $$\sum_{ l =1}^{n_k}e_{il}^k(L(g))e_{lj}^k(L(h))= \sum_{ l =1}^{n_k}e_{li}^k(g)e_{jl}^k(h).$$
 
Then 
$$( \rho_g \circ L^\ast )(e_{ij}^k)= \sum_{ l =1}^{n_k}e_{li}^k(g)e_{jl}^k.$$

Since

$$ Tr_k(\rho_g \circ L^\ast )= \sum_{1 \leq i,j \leq n_k} \langle (\rho_g \circ L^\ast )(e_{ij}^k), e_{ij}^k  \rangle= \sum_{1 \leq i,j \leq n_k}\sum_{ l =1}^{n_k}e_{li}^k(g) \langle e_{jl}^k,e_{ij}^k  \rangle$$

\noindent 
and
$\langle e_{jl}^k,e_{ij}^k \rangle = \delta_{(j,l), (i,j)}$,
%is equal to $ 1$ when $ j=i,l=j $ and is equal to $0 $ otherwise, 
we get 

$$ Tr_k(\rho_g \circ L^\ast )= \sum_{i=1}^{n_k}e_{ii}^k(g)=\chi_k(g) $$

Let us suppose now that $k$ is such that $ \varepsilon_k = -1. $ Let $ n_k = 2m_k; $ then we
can find a symplectic basis $\mathcal{U}_k^{-}= \{ u_i, 1 \leq i \leq n_k \} $
 of $ U_k $ such that $\langle u_i,u_{i+m_k} \rangle =1 $ 
 and $\langle u_i,u_{j+m_k} \rangle =\langle u_{i+m_k},u_{j+m_k} \rangle=
\langle u_i,u_j\rangle = 0,    i  \neq j  \;\; ( 1 \leq  i,j \leq m_k). $

Equation  1 gives us now the following relations for the
matrix coefficients $ e_{ij}^k(g) $ of $ \pi_k(g) $
 with respect to this basis:
$$\langle \pi_k (g) (u_i), u_j\rangle= -e_{j+m_k \;i}^k (g), \; \; j \leq m_k , $$
$$\langle \pi_k (g) (u_i), u_j\rangle=e_{j-m_k\; i}^k (g), \; \;  j > m_k$$

 Therefore
 
 $$\langle \pi_k (L(g)) (u_i), u_j\rangle =  - \langle \pi_k (g) (u_j), u_i\rangle=e_{i+m_k
 \; j}^k(g), \; \; i \leq m_k , $$ 
 and
 
 $$\langle \pi_k (L(g)) (u_i), u_j\rangle= -e_{i-m_k\; j}^k(g), \; \; i > m_k$$.
 
 Taking into account these relations and equation 2 we get the following relations between the matrix coefficients of $\tilde{L} (\pi_k(g)) = \pi_k(L(g)) $ and $ \pi_k(g):$
  
\begin{equation}
e_{ij}^k(L(g)) =e_{j+m_k \; i-m_k}^k(g), \; \; j \leq m_k , i > m_k
\end{equation}

\begin{equation}
e_{ij}^k(L(g)) = -e_{j-m_k 
\;  i-m_k}^k(g), \; \; j > m_k  , i > m_k
\end{equation}
\begin{equation}
e_{ij}^k(L(g)) = -e_{j+m_k 
\;  i+m_k}^k(g), \; \; j \leq m_k  , i \leq m_k
\end{equation}
\begin{equation}
e_{ij}^k(L(g)) = e_{j-m_k \;  i+m_k}^k(g), \; \; j > m_k  , i \leq m_k
\end{equation}

So we get:

For  $i,j \leq m_k$
\begin{equation}
(\rho_g \circ L^\ast )(e_{ij}^k)= \sum_{ l =1}^{ m_k }-e_{l+m_k
\;  i+m_k}^k(g)e_{j+m_k
\; l+m_k}^k + \sum_{l=m_k+1}^ { n_k }e_{l-m_k 
\; i+m_k}^k (g)e_{j+m_k \; l-m_k}^k.
\end{equation}
For $ i \leq m_k, j > m_k $
\begin{equation}
(\rho_g \circ L^\ast )(e_{ij}^k )= \sum_{ l=1}^ { m_k }-e_{l+m_k 
\;  i+m_k}^k (g)e_{j-m_k \; l+m_k}^k + \sum_{l=m_k+1}^ {n_k }e_{l-m_k \; i+m_k}^k (g)(-e_{j-m_k \; l-m_k}^k ).
\end{equation}
For $ i > m_k, j \leq m_k $
\begin{equation}
(\rho_g \circ L^\ast )(e_{ij}^k )= \sum_{ l =1}^{ m_k }e_{l+m_k \; i-m_k}^k (g)(-e_{j+m_k 
\;  l+m_k}^k) + \sum_{l=m_k+1}^ {n_k }-e_{l-m_k \; i-m_k}^k)(g)e_{j+m_k \; l-m_k}^k.
\end{equation}
And for $ i > m_k, j > m_k $
\begin{equation}
(\rho_g \circ L^\ast )(e_{ij}^k) = \sum_{ l=1} ^{m_k }e_{l+m_k
\;  i-m_k}^k (g)e_{j-m_k \; l+m_k}^k  + \sum_{l=m_k+1}^ { n_k }-e_{l-m_k \; i-m_k}^k (g)(-e_{j-m_k \; l-m_k}^k ).  
\end{equation}

Notice that:
\begin{itemize}
\item[a)]
If $ i,j \leq m_k $ then $e_{j+m_k \; l+m_k}^k \neq e_{ij}^k $ and $e_{j+m_k \; l-m_k}^k \neq e_{ij}^k $
\item[b)]
If $ i \leq m_k ,\;  j > m_k $ then $e_{j-m_k
\; l+m_k}^k =e_{ij}^k $ if and only
if \mbox{$ j=m_k+i>m_k $} and $ l+m_k= j>m_k, $   moreover $ e_{j-m_k \; l-m_k}^k \neq e_{ij}^k $
\item[c)]
If $ i > m_k , \; j \leq m_k $ then $e_{j+m_k \; l-m_k}^k = e_{ij}^k $ if and only if $
j+m_k=i>m_k $ and $ l-m_k = j<m_k, $ moreover  $ e_{j+m_k \; l+m_k}^k \neq e_{ij}^k $
\item[d)]
If $ i > m_k , \; j > m_k $ then $e_{j-m_k \; l+m_k}^k  \neq e_{ij}^k $ 
 and $ e_{j+m_k \; l-m_k}^k \neq e_{ij}^k $
\end{itemize}

Therefore
$$ Tr_k(\rho_g \circ L^\ast )= \sum_{j=m_k+1}^ {n_k}-e_{jj}^k (g) + \sum_{j=1} ^{ m_k}-e_{jj}^k (g) = -\chi_k(g).
$$

Hence, taking into account both cases, we get

$$ Tr(\rho_g \circ L^\ast )=\sum_{k=1}^{r} Tr_k(\rho_g \circ L^\ast )= \sum_{k=1}^{r}  \varepsilon_k \chi_k(g),
$$
which proves iv.
 \end{proof}
 
\begin{defi}
For any mapping $L$ from $G$ to itself, we denote by  $Fix_G(L)$ the fixed point set   $ \{g \in G : L(g) = g \}$ 
of  $L.$

\end{defi}

\begin{theo}
If $L$ is an involutive antiautomorphism of $G$ such that  
$   L(g)$ is conjugated to $g$ for all $g \in G$  and
$$ |Fix_G(L)|= d_G. $$
%\sum_{k=1}^ { r } n_k,  $$ 

\noindent
then
$$ Tr(\rho_? \circ L^\ast )= \chi_G $$ 
%\sum_{k=1}^ {r}   \chi_k(g)
  and $$L^\ast  = T. $$

\end{theo}
%%%%%.   CHanges begin here ....
%.  ask Antonio about the wild conjecture :   A4 is the only non commutative finite group to have a natural Gelfand model  !!!!

\begin{proof} From Theorem 2.iv) we have: $ Tr(\rho_g \circ L^\ast ) = \Sigma_{k=1}^ { r} \varepsilon_k \chi_k (g). $

 If  we evaluate $ Tr(\rho_g \circ L^\ast ) $ on $ e$ we
obtain
 $$ Tr ( L^{\ast }) = \sum_{k=1}^ { r} \varepsilon_k n_k $$

Furthermore,  from Theorem 2.iii) 
$$ Tr(\rho_e \circ L^\ast ) = |\{h \in G : h^{-1}L(h)) = e \}|= |\{h \in G : L(h)) = h \}|  =|Fix_G(L)| $$

\noindent
but $   Tr(\rho_e \circ L^\ast )= Tr ( L^{\ast }) $ 
and by hypothesis 
$$ |Fix_G(L)|= d_G = \sum_{ k =1}^{r } n_k. $$ 

 So, \;  $  \sum_{k=1}^ { r} \varepsilon_k n_k = \sum_{ k =1}^{
r } n_k; $
since $ n_k > 0 $, we conclude that $\varepsilon_k =
1$ for all $ k,$ and $ Tr(\rho_g \circ L^\ast )= \chi_G(g) \; (g\in G). $

 Then, using equation 3 it follows  that 
$$ L^{\ast }(e_{ij}^k)(g) = e_{ij}^k(L(g)) = e_{ji}^k(g), $$  
hence $L^\ast  = T. $
    
\end{proof}
%prop 3. %%%
\begin{prop}

If $L$ is an involutive antiautomorphism of $G$ such that  $  L(g) $ is conjugate to $g$ for all   $g \in G, $ 
  then   
\begin{equation}
  Tr(\rho_g \circ L^\ast )=\sum_{1 \leq k \leq r} (\frac{1}{|G|} \sum_{h \in G }\chi_k(L(h)h^{-1}) )  \chi_k(g)
\end{equation}
\end{prop}

\begin{proof}

To prove  equation (12) we compute the Fourier coefficients $
\lambda_k , 1 \leq k \leq r $ of the central function $$Tr(\rho_g
\circ L^\ast ) = Tr( L^\ast  \circ \sigma_{L(g)^{-1}}) $$ with
respect to the basis $$ \{\chi_k =  \sum_{ i=1} ^{n_k } e_{ii}^k
:1 \leq k \leq r \}$$  of the centre $Z$ of  $ \mathbb{C}[G].$ 
First we observe that

$$( L^\ast  \circ \sigma_{L(g)^{-1}})(e_{ij}^k)(h) = e_{ij}^k(L(g)L(h)) =\sum_{l=1 }^
 {n_k } e_{il}^k(L(g))( e_{lj}^k \circ L )(h)$$
 
\noindent 
and

$$ \langle (L^\ast  \circ \sigma_{L(g)^{-1}})(e_{ij}^k), e_{ij}^k
\rangle = \frac{1}{|G|} \sum_{h \in G }\sum_{ l=1}^ {n_k }
e_{il}^k(L(g)) e_{il}^k (L(h)) e_{ji}^k ( h^{-1}). $$

therefore

$$Tr( L^\ast  \circ \sigma_{L(g)^{-1}})  =\sum_{k=1}^{r}
\sum_{1\leq i, j \leq n_k }\frac{1}{|G|} \sum_{h \in G }\sum_{l=1}^{ n_k } e_{il}^k(L(g)) e_{lj}^k (L(h)) e_{ji}^k ( h^{-1})
$$

Since
 $$ \sum_{ j=1} ^{n_k}e_{lj
}^k(L(h)) e_{ji}^k(h^{-1})=
e_{li}^k (L(h))h^{-1}) $$
 
we get that

$$Tr( L^\ast  \circ \sigma_{L(g)^{-1}})  = \sum_{ k =1}^{r}\sum_{1 \leq i,l \leq n_k }\frac{1}{|G|} \sum_{h \in G } 
 e_{li}^k(L(h)h^{-1}) (e_{il}^k \circ L)(g)$$
 
and then, for    $ 1 \leq k' \leq r,  $

$$
\lambda_{k'} = \langle\sum_{ k=1}^ {r}\sum_{1 \leq i,l \leq n_k
}\frac{1}{|G|} \sum_{h \in G } 
 e_{li}^k(L(h)h^{-1}) (e_{il}^k \circ L), \chi_{k'} \rangle
$$

By hypothesis $ \chi_k(L(g))= \chi_k(g), 1 \leq k \leq r  $, so

$$
\lambda_{k'} = \sum_{ k=1} ^{ r}\sum_{1 \leq i,l \leq n_k
}\frac{1}{|G|} \sum_{h \in G } 
 e_{li}^k(L(h)h^{-1}) \sum_{ s=1} ^{n_{k'}}(\frac{1}{|G|} \sum_{g \in G } e_{il}^k (L(g))\overline{e_{ss}^{k'} (L(g))},
$$

But $\frac{1}{|G|} \sum_{g \in G } e_{il}^k (L(g))\overline{e_{ss}^{k'}
(L(g))}= \langle e_{il}^k ,e_{ss}^{k'}\rangle, $  then

$$
\lambda_{k'} = \sum_{ s=1}^ {n_{k'}}\frac{1}{|G|} \sum_{h \in G }
 e_{ss}^{k' }(L(h)h^{-1})= \frac{1}{|G|} \sum_{h \in G }
 \chi_{k'}(L(h)h^{-1}).
$$
\end{proof}

\begin{prop}
 
Let  $L$ be an involutive antiautomorphism of $G$ such that $ L(g) $ is conjugate to  $g$ for all $ g \in G. $ Then the
following conditions are equivalent
\begin{itemize}
\item[i.]
$$  |Fix_G(L)|= d_G $$
\item[ii.]
$$  \frac{1}{|G|}\sum_{h \in G }\chi_k(L(h)h^{-1})= 1, \ \ (1 \leq k \leq r )$$
\end{itemize}

\end{prop}

\begin{proof}
It follows from theorem 2 and proposition 2.
 
\end{proof}

\begin{prop}

Let $L$  an involutive antiautomorphism of $G$ ,  $ \tau : G \rightarrow G $ defined by $
\tau(g) = L(g^{-1}) $ and let $ c_\tau (\chi) = |G|^{-1}\sum _{g \in G} \chi ( g \tau(g)) $ be the twisted Frobenius - Schur indicator defined by Kawanaka and Matsuyama \cite{K-M}. 

Then 
$c_\tau (\chi_k) =  \frac{1}{|G|} \sum_{h \in G } \chi_{k}(hL(h^{-1})
=\pm 1 $ if 
$\chi_k(L(g))= \chi_k (g). $ 

Moreover if   the matrix representation  $
R_k $ corresponding to  $\chi_k$  satisfies $ R_k(L(g))= R_k(g), $  then $
c_\tau(\chi_k) = 1.$

\end{prop}
\begin{proof}
 Due to $ g \tau(g)) = g L(g^{-1}),$  we have

$$ \sum_{ g \in G}\chi(g L(g^{-1}))= \sum_{ g \in G}\overline{\chi(L(g)
g^{-1})} $$

and  since $$Tr(\rho_g \circ L^\ast )= \overline{Tr(\rho_g \circ
L^\ast )}$$

we can write:

$$Tr(\rho_{g^{-1}} \circ L^\ast )= \sum_{ k=1}^ { r}\frac{1}{|G|} \sum_{h \in G }
\overline{\chi_{k}(L(h)h^{-1})\chi_k(g^{-1})}, $$

 and therefore

$$Tr(\rho_{g} \circ L^\ast )= \sum_{ k=1}^ { r}\frac{1}{|G|} \sum_{h \in G }
\chi_{k}(hL(h^{-1}))\chi_k(g). $$

Our proposition follows. 

\end{proof}

%Notice that condition ii. above does not follow from condition i (Lg being conjugated to g)  .   Indeed,  for $G = SL(2,\mathbb{F}_q) $ with $q$ congruent to  $1 \;\;  mod \;4 $  we have that  i. holds for   $L(g) = g^{-1} \;\;\; (g \in G)$ but  $L$ has only two fixed points in $G.$ 

\section{  Applying the twisted trace to test  natural Gelfand Model  candidates}

%%%. test wether the Gelfand character could be a permutation character  associated to a naturally suggested X. %% 

% Remark for G = SL2, q the arithmetic obstruction vanishes.      dim M = q(q+1) |  (q-1)q(q+1)  =. order of G   holds 
%but   G has non trivial centre and when you average the irred characters over the inducing subgroup, which is the diagonal torus, you find as multiplicity the summation of the irred char over the centre, which is 0 or 2 instead of 1, according to the char parameter being a square or not.So we are done if the centre is trivial.

% Conjecture    dim M divides order of the group  G.  and trivial centre implies % G has a natural Gelfand Model. 

Let $ G $  be a group for which there exists  an involutive  antiautomorphism $L$  satisfying all the  conditions of theorem 3 ,
%%. where not being able to tell L(g) from g. with the help of the characters is equivalent to they being conjugate 
Then we can consider the  $ G $ - set  $   X  =   Fix_G(L) $
where the action of $G$ in $X$ is given by        $ g \cdot h = g^{-1}hL(g)^{-1}. $

Let us call $\tau$ the associated natural representation of  $G$ in $L^2(X)$.

Let  $ Z  $ denote the centre  $Z(G)$ of  $G$;    we note that  if $ z \in Z$ then  $ z \in X $ because  $L(z) $ is conjugated to $z$ by hypothesis. Furthermore  since $ \chi_ \tau ( g) = | \{h \in X : ghL(g)=h \}|    $
 for all   $ g \in  G $,  it follows that   $     \chi_ \tau (z)  =  |  X |   $   for  all   {\em involutions}
 $ z \in Z $.  So the character $\chi_ \tau$ does not tell central involutions from the identity.

Since the dimension of $ L^2(X)$ is equal to the Gelfand dimension $d_G$ of $G$, the representation $ (L^2(X), \tau) $  {\em could} be a Gelfand Model for $ G $. 
To check  whether  this is indeed the case,   we can of course compare its character $  \chi_ \tau $ with the Gelfand character $\chi_G$ that we have described as a twisted trace: $\chi_{G}(g) =Tr( \rho_g \circ L^*) , \;g \in G $.

We see then that

\begin{prop}
Let $ G $  be a group for which there exists  an involutive  antiautomorphism $L$  satisfying all the  conditions of theorem 3 and consider the  $ G $ - set  $   X  =   Fix_G(L) , $ where the   action of $G$  is given by   $ \;\; g \cdot h = g^{-1}hL(g)^{-1} \;\;\; (g \in G, h \in X). $ 

Then a necessary condition for  $(L^2(X), \tau) $ to be a Gelfand Model is  that
 
%$$  \chi_ G(e) =  \chi_ \tau (e)  \; ( = | X |  ) $$  

 $$  | X | = d_G $$
and
     $$     \chi_ G(z) =    |  X  |  $$  
 for all central involutions   $z$, or equivalently,  in terms of $L$,    
$$      |  \{ g \in G  :  g^{-1}L(g) = z   \}| =  |  X    |   $$ .

\end{prop}

\subsection{  The case of $ G =  SL(2, \mathbb{F}_q),\;$  $q$ odd}

As an illustrative example we consider here the case of  $ G =  SL(2, \mathbb{F}_q),\;$  $q$ odd.

  Recall that the Gelfand dimension of $G$, i. e. the sum of the dimensions of all its irreducible representations, is $ q(q+1)$.

 For   $  k =   \left( \begin{smallmatrix} 1 & 0 \\
0 & -1 \end{smallmatrix} \right),$   we may define   $ L: G\rightarrow G, L(g)= kg^{-1}k $.

 We note that if $ g= \left(\begin{smallmatrix} a & b \\ c & d \end{smallmatrix} \right) \in G $ then 
$ L(g)=  \left(\begin{smallmatrix} d & b \\ c & a \end{smallmatrix} \right), $
 and $$ X=Fix_G(L)= \{\left(\begin{smallmatrix} a & b \\ c & a \end{smallmatrix} \right)\in G : a^2-bc=1\}.$$.

\begin{prop}  $ L $ is an involutive antiautomorphism of $G$ satisfying:
\begin{itemize} 
\item[i.] $ L(g)$ is conjugated with $ g $, for all   $g \in G.$
\item[ii.] $|X|= q(q+1)$  
\item[iii.] The action of $ G $ on   $ X $ is transitive.
\end{itemize}
\end{prop}
\begin{proof}     
i) This is clear since $L(g) $ has same trace and determinant as  $g$ and it is a scalar matrix if and only if $g$ ia, for all  $ g \in G.$

ii) For   $   \left(\begin{smallmatrix} a & b \\  c & a  \end{smallmatrix}
\right)  \in   X, $    we must have   $bc= a^2 - 1$.  Now,   if  $a^2 - 1 = 0$, i. e.  $ a = \pm 1,$   we  have  $2q-1$  pairs  $(b,c)$ such that  $bc =  a^2 - 1.$  On the other hand, if    $a^2 - 1 \neq  0, $ then  there are $ q-1 $  pairs  $(b,c)$ such that $bc  = a^2 - 1.$     
So 
$$   |X|  =   2 \cdot (2q-1) +  (q-2) \cdot  (q-1)  =  q(q+1). $$ 

iii)We have that $ | Orb_G (e)|= {|G|} / {|Stab_G (e)|} $ and $ Stab_G (e)= \{g \in G :L(g) = g^{-1} \}. $ 
 For  $ g =  \left(\begin{smallmatrix} a & b \\
c & d
\end{smallmatrix}
\right) \in G, $ this means     $ b=c=0 $ and $ d=a^{-1}$ Therefore $ |Stab_G (e)| = | \{ \left(\begin{smallmatrix} a & 0 \\
0 & a^{-1}
\end{smallmatrix}
\right):a \in \mathbb{F}_q^* \}|    = q-1. $ 

So  $ | Orb_G (e)|= \frac{q(q^2-1)}{q-1}= q^2 +q = |X|$  and  therefore the action is transitive.
\end{proof}

In this way we get that $ L $ is an involutive antiautomorphism of $G= SL(2,\mathbb{F}_q) $ satisfying all  conditions in Theorem  1   above.   Hence  for the central involution $z= \left(\begin{smallmatrix} -1 & 0 \\
0 & -1
\end{smallmatrix}
\right),$ we have $$ \chi_ G(z) = |\{ h \in G: h^{-1}L(h)= z \}|.$$ 
But if $ h=\left(\begin{smallmatrix} a & b \\
c& d
\end{smallmatrix}
\right) \in G$ satisfies $ h^{-1}L(h)= z$ then  $$ \left(\begin{smallmatrix} d & b \\
c & a
\end{smallmatrix}
\right)= \left(\begin{smallmatrix} -a & -b \\
-c & -d
\end{smallmatrix}
\right) $$ 
 and that is equivalent to  $b=c=0 $ and $ -a^2=1.$

So,  if $ q \equiv 3 \;\;  mod (4)$ then $\chi_ G(z) = 0 $ and  if $ q \equiv 1 \;\;  mod (4)$ then $\chi_ G(z) =  2 $.
 In both cases $      \chi_ G(z) \neq  \chi_ \tau (z).   $ Therefore  $ (L^2(X), \tau) $ is not a Gelfand Model for $G$.
 
% \end{proof}

 Motivated by the previous example, we may notice the following.

 \begin{lema}
 
  For every  central element  $z \in G$ of order $2,$ we have. 
 $$  \chi_G(z ) \neq  \chi_G(e).   $$
 
 \end{lema}

\begin{proof}
Indeed,   notice that    $\chi_\pi (z)  =   \pm \; dim(\pi )   $   for all   irreducible representations     $\pi$  of $G$    
 and so for    $  \chi_G(z  ) = \chi_G(e)    $ to hold  we need that  
  $\chi_\pi (z)  =  dim(\pi )    =  \chi_\pi (e)$
    for all   irreducible representations $\pi,$   
 but this is absurd since   the regular character  $\chi_\rho $ tells    $z$ from   $e,$ because   
  $\chi_\rho (e)  = |G|   \neq   0  =  \chi_\rho (z).   $
 
 \end{proof} 
 
 \begin{prop}  Let $G$ a finite group  and an  involutive anti automorphism   $L$ of G such that  
 
 \begin{itemize}
\item[i.]
$L(g)$  is conjugated to  $  g, $ for all $g \in G;$

\item[ii.]
$|Fix_G(L)|= d_G$

\end{itemize}
Then  
  if $G$ has central (non trivial) involutions we conclude that   the natural representation $ (L^2(X), \tau) $ associated to    the $G-$ space $X = Fix_G(L)$  cannot be a Gelfand Model for $G.$       
 \end{prop}
 \begin{proof}
 This follows immediately from the necessity of the condition   
   $     \chi_ G(z) =  \chi_ \tau (z)  $  
 for all central involutions   $z$  and lemma 1.
 \end{proof} 

%\section {Non transitivity of the   action of $G$ on $ Fix_G(L)  $ as an obstruction to the existence of an involution fixed point geometrical Gelfand Model }

\subsection{  The case of the quaternion group $  Q$}

The quaternion group $Q$ is an  illustrative example of a group admiting   involutive  antiautomorphisms $L$  satisfying all the  conditions of theorem 3 whose  fixed point set 
  $  X =  Fix_Q(L)  $, endowed with the natural twisted action of $Q$,  does not afford a Gelfand model for $Q$, because  $Q$ has non trivial central involutions, but such that the $Q$-set $X$ is {\em non transitive}.  
  
Recall that the quaternion group $Q$ may be described as  
$$ Q = \{1, -1, i, -i, j, -j, k, -k\} $$
where $g^2 = -1$ for all $g \in Q,  g \neq \pm 1$, and 
$$ ij = k = -ji, jk = i = -kj, ki = j = - ik $$ 

% $$ \mathbb{Q}= \langle i\;\;,j\;\;| i^4=e,\;\; i^2=j^2,\;\;jij^{-1}=i^{-1}\rangle$$

  Every element is only conjugated to its inverse and the centre of $G$ is $ \{1, -1\}. $ The five conjugacy class of $G$ are: 
  %$ [e], [i], [j],[ij], [i^2] $.   
  $ \{1\}, \{ -1   \}, \{ \pm i  \}, \{ \pm j  \},
  \{ \pm k  \}.$
  
The quaternion group $Q$ has four 1-dimensional irreducible representations and one 2-dimensional  irreducible representation. So  $d_Q = 6$.
  
The posibilities for $ L(i)$ and  $L(j) $  for $L$ to be  an  involutive {\em antiautomorphism}   of $G$ such that $ L(g) $ is conjugated whith $ g, $ for all $g \in Q $ are:

$ L(i) = \pm i $   and   $ L(j)= \pm j $.
 
  If we define $ L(i)=i$ and $L(j)=j$ then  $L(k) = - k $  and 
  $$ Fix_Q(L)  = \{\pm 1, \pm i, \pm j \}  $$    
  
  If we define $ L(i)=i$ and $L(j)= -j $   then  $ L(k) = k $ and 
 
 $$ Fix_G(L)  = \{\pm 1, \pm i, \pm k \}  $$   
 
 Analogously, if we define $ L(i)= -i$ and $L(j)= j $   then  $ L(k) = k $ and
 
 $$ Fix_Q(L)  = \{\pm 1, \pm j, \pm k \}  $$

 So, in all these cases we get  $| Fix_Q(L)| = 6 = d_Q.$ 
 
Notice that if we define   $ L(i)= -i$ and $L(j)= -j, $   then  $ L(k) = -k $ and 
 
 $$ Fix_Q(L)  = \{\pm 1 \}  $$

 Notice that any of the first three  antiautomorphisms defined above  clearly  satisfies all the conditions of Theorem 3. We have then for  the Gelfand Character of $Q$:
  $$ \chi_ Q(g) = |\{ h \in Q: h^{-1}L(h)= g \}|.$$
  
 It follows that for the non trivial central involution  $-1$ we have, for the first such $L$ (the other two cases being obtained by cyclic permutation of $i,j ,k),$: 
  $$ \chi_ Q(-1) = |\{ h \in G: L(h)= -h \}| = |\{\pm k  \}| = 2 \neq 6 = d_Q $$
  So $ \chi_ Q $ is different from the character $\chi_ \tau $ of $Q$ associated to its natural 
  action on $ Fix_Q(L),$ given by $g.x = g^{-1}xL(g)^{-1}.$
  
  Notice however that the   action of $Q $ on   $X=Fix_Q(L)$ is not transitive. There are in fact $ 3 $ orbits of $Q$ in $X$, to wit :

 $$  Orb_Q (1)=\{ \pm 1 \}, \; Orb_Q (i)=\{ \pm i \}, \; Orb_Q (j)=\{ \pm j \}.$$

\bigskip 
 {\bf \large Acknowledgements.}

 The authors thank Pierre Cartier for helpful discussions related to  this work.

\vspace{2cm}

%\noindent
%{\em Acknowledgements}:    

\end{document}